\documentclass[a4paper,11pt]{article}
\usepackage{latexsym}
\usepackage{amssymb}
\usepackage{enumerate}
\usepackage{amsfonts}
\usepackage{amsmath}
\usepackage[dvipdfmx]{graphicx} 
\usepackage[paper=a4paper,left=30mm,right=20mm,top=25mm,bottom=30mm]{geometry}

\newcommand{\qed}{$\Box$}

\newenvironment{@abssec}[1]{%
    \if@twocolumn

      \section*{#1}%
    \else

      \vspace{.05in}\footnotesize
      \parindent .2in
 {\upshape\bfseries #1. }\ignorespaces
    \fi}

    {\if@twocolumn\else\par\vspace{.1in}\fi}

\newenvironment{keywords}{\begin{@abssec}{\keywordsname}}{\end{@abssec}}

\newenvironment{AMS}{\begin{@abssec}{\AMSname}}{\end{@abssec}}

\newcommand\keywordsname{Key words}
\newcommand\AMSname{AMS subject classifications}
\newcommand\AMname{AMS subject classification}
\newtheorem{theorem}{Theorem}
 \newtheorem{lemma}[theorem]{Lemma}
 
 \newtheorem{proposition}[theorem]{Proposition}

\def\qed{\vbox{\hrule height0.6pt\hbox{%
  \vrule height1.3ex width0.6pt\hskip0.8ex
  \vrule width0.6pt}\hrule height0.6pt
 }}

\title{Large time behavior of temperature in two-phase heat conductors\thanks{This research was partially supported by the Grant-in-Aid
for Scientific Research (B) ($\sharp$ 18H01126)  of
Japan Society for the Promotion of Science and NRF (of S. Korea) grants No. 2019R1A2B5B01069967.}}

\author{Hyeonbae Kang\thanks{Department of Mathematics and Institute of Applied Mathematics, Inha University, Incheon 22212, S. Korea (hbkang@inha.ac.kr).}  \and Shigeru Sakaguchi\thanks{Research Center for Pure and Applied Mathematics, Graduate School of Information Sciences, Tohoku University, Sendai, 980-8579, Japan (sigersak@tohoku.ac.jp).}}

\date{}
\begin{document}
\maketitle

\begin{abstract}
We consider the Cauchy problem for the heat diffusion equation in the whole Euclidean space consisting of two media with different constant conductivities. The large time behavior of temperature, the solution of the problem, is studied when initially temperature is assigned to be 0 on one medium and $1$ on the other. We show that under a certain geometric condition of the configuration of the media, temperature is stabilized to a constant as time tends to infinity. We also show by examples that temperature in general oscillates and is not stabilized.
  \end{abstract}



\pagestyle{plain}
\thispagestyle{plain}


\section{Introduction}\label{introduction}

This paper concerns  the Cauchy problem for the heat diffusion equation in the whole Euclidean space which is occupied by two heat conducting media with different constant conductivities. It deals with the question of stabilization of temperature (the solution of the problem) as time tends to infinity when the initial data is given by the characteristic function of a medium.  The reason why we introduce such a specific initial data is to clarify geometry of the composite media. We find a condition on the configuration of the media under which temperature converges to a constant as time tends to infinity. We also show by examples that temperature in general oscillates and does not converge.

The stabilization problem of the solution of linear parabolic equations as time tends to infinity has been studied by many authors. For example, Kamin considered in \cite{K1976} the Cauchy problem for the uniformly parabolic diffusion equations with bounded initial data. Under the assumption that the diffusion coefficients converge to the Kronecker symbol at infinity (the condition (9) in \cite{K1976}), that is, the equation converges to the exact heat equation  at infinity, it is shown that the solution converges to a constant as time tends to  infinity if and only if the averages of the initial data over balls converge to a constant as the radii of the balls tend to infinity. A generalization of this result to some degenerate parabolic equation was given in \cite[Theorem 1.1]{EKT2009}. We refer to \cite{E1969, Z1977} and references therein for further work on stabilization of the parabolic equations. We also mention that the oscillatory behavior of solutions of the Cauchy problem for the exact heat equation with unbounded initial data was studied in \cite{CT2020}.

To present the results of this paper precisely, let $\Omega$ be a domain in $\mathbb R^N$ with $N \ge 2$ so that $\Omega$ and $\mathbb R^N \setminus \Omega$ constitute the two media. We suppose that $\partial\Omega\not=\emptyset$ and $\partial\Omega$ is connected.
Denote by $\sigma=\sigma(x)\ (x \in \mathbb R^N)$ the heat conductivity distribution of the whole medium given by
\begin{equation}
\label{conductivity constants}
\sigma =
\begin{cases}
\sigma_+ \quad&\mbox{in } \Omega, \\
\sigma_- \quad &\mbox{in } \mathbb R^N \setminus \Omega,
\end{cases}
\end{equation}
where $\sigma_-, \sigma_+$ are positive constants such that $\sigma_-  \not= \sigma_+$. We consider the Cauchy problem for the heat diffusion equation:
\begin{equation}\label{heat Cauchy}
\begin{cases}
u_t =\mbox{ div}(\sigma \nabla u)\quad &\mbox{in }  \mathbb R^N\times (0,+\infty), \\
u\ ={\mathcal X}_{\Omega}\ &\mbox{on } \mathbb R^N\times \{0\},
\end{cases}
\end{equation}
where ${\mathcal X}_{\Omega}$ denotes the characteristic function of the set $\Omega$. We look for a geometric condition on $\Omega$ such that the unique bounded solution  $u =u(x,t)$ to \eqref{heat Cauchy} tends to a constant as $t \to \infty$.
As far as we are aware of, the question of stabilization for the multi-phase heat conductors has not been considered before. Recently, a geometric question related to the diffusion over such multi-phase heat conductors has been dealt with in \cite{CSU2019, S2020}.

To gain a better understanding of the condition for the stabilization given later in \eqref{cone}, let us first consider conic regions. Let $S^{N-1}$ be the unit sphere in $\mathbb R^N$. For a domain $A$ of $S^{N-1}$ such that $\partial A \not=\emptyset$, we set $\Omega_A$ to be the cone over $A$, namely,
\begin{equation}
\label{definition of cone}
\Omega_A = \{ x \in \mathbb R^N : x = r\omega,\ r > 0,\ \omega \in A \}.
\end{equation}
Denote by $\sigma_A=\sigma_A(x)\ (x \in \mathbb R^N)$  the conductivity distribution of the whole medium given by
\begin{equation}
\label{conductivity constants selfsimilar}
\sigma_A =
\begin{cases}
\sigma_+ \quad&\mbox{in } \Omega_A, \\
\sigma_- \quad &\mbox{in } \mathbb R^N \setminus \Omega_A.
\end{cases}
\end{equation}

The following proposition can be proved using the self-similarity of the solution (see the proof in the next section).

\begin{proposition}\label{thm:main1}
Let $u_A =u_A(x,t)$ be the unique bounded solution of problem \eqref{heat Cauchy} with $\Omega, \sigma$ replaced by $\Omega_A, \sigma_A$, respectively. It holds that $0< u_A(0,1) <1$ and
\begin{equation}
\label{stabilization for self-similar solutions}
\lim\limits_{t \to \infty} u_A(x,t) = u_A(0,1)
\end{equation}
uniformly in $x$ belonging to any fixed compact set in $\mathbb R^N$.
\end{proposition}

We now present a condition on the shape of the domain which is sufficient for stabilization of the solution.  Let $A$ be a domain of $S^{N-1}$ such that there is a point $p\in A$ with $-p\not\in A$. For $A \subset S^{N-1}$ and such a point $p$, we say $A$ is {\it starshaped} with respect to $p \in A$  if  for every point $\omega \in A$ the shortest geodesic connecting $\omega$ and $p$ in $S^{N-1}$ is contained in $A$, or equivalently, if $t\omega + sp \in A$ for any $\omega \in A$ and a pair of nonnegative numbers $t,s$ such that $t\omega + sp \in S^{N-1}$. If $A \subset S^{N-1}$ is  starshaped with respect to $p \in A$, then
\begin{equation}
\label{cone2}
\Omega_A \subset  \Omega_A - s p
\end{equation}
for all $s > 0$. Here and throughout this paper $\Omega_A - y$ denotes the translate of $\Omega_A$ by $y$, that is,
$$
\Omega_A - y = \{ x-y  : x \in \Omega_A \}.
$$
In fact, if $x \in \Omega_A$, then $x=t\omega$ for some $t>0$. Since $A$ is starshaped with respect to $p$, $x+sp \in \Omega_A$ for all $s>0$, and hence \eqref{cone2} follows.

The sufficient condition for the stabilization of the solution to \eqref{heat Cauchy} to hold on the domain $\Omega$ is as follows:
\begin{equation}
\label{cone}
\Omega_A \subset \Omega \subset  \Omega_A - h p
\end{equation}
for some $h > 0$. Roughly speaking, this condition means $\partial \Omega$ lies in $( \Omega_A - h p) \setminus \Omega$, but is of arbitrary shape. For such domains we have the following theorem.

\begin{theorem}\label{stabilization}
If the domain $\Omega$ in $\mathbb R^N$ satisfies \eqref{cone} for some $h>0$,  $A \subset S^{N-1}$ and  $p \in A$ with $-p\not\in A$, where $A$ is starshaped with respect to  $p$,  then  the unique bounded solution $u$ of problem \eqref{heat Cauchy} satisfies that
\begin{equation}\label{limita}
\lim\limits_{t \to \infty} u(x,t) = u_A(0,1)
\end{equation}
uniformly in $x$ belonging to each compact set in $\mathbb R^N$.
\end{theorem}

The second theorem gives an example of oscillatory behavior of the solutions of  problem \eqref{heat Cauchy}. For that we deal with general conductivity distributions $\sigma=\sigma(x)$, not necessarily two-phase conductivity distribution.
\begin{theorem}\label{oscillation}
Let $m \le M$ be positive constants. There exists a domain $\Omega$ in $\mathbb R^N$ such that for any conductivity $\sigma=\sigma(x)\ (x \in \mathbb R^N)$ satisfying
\begin{equation}\label{mM}
0 < m \le \sigma(x) \le M\ \mbox{ for every } x \in \mathbb R^N,
\end{equation}
the unique bounded solution $u$ of problem \eqref{heat Cauchy} satisfies that
\begin{equation}
\label{oscillatory behavior}
0< \liminf\limits_{t \to \infty} u(0,t) < \limsup\limits_{t \to \infty} u(0,t) < 1.
\end{equation}
\end{theorem}

We remark that since Theorem \ref{oscillation} is proved only using the Gaussian bounds for the fundamental solutions of diffusion equations (see \eqref{Gaussian bounds}), it holds also for the diffusion equations of the form
$$
u_t = \sum_{i,j=1}^N\frac \partial{\partial x_i}\left(a_{ij}(x,t)\frac{\partial u}{\partial x_j}\right)
$$
where the coefficients satisfy the following for some positive constants $m, M$
$$
m |\xi|^2\le\sum_{i,j=1}^Na_{ij}(x,t)\xi_i\xi_j\le M |\xi|^2\ \mbox{ for } x,\xi\in \mathbb R^N\mbox{ and } t>0.
$$
Such diffusion equations have been dealt with in \cite{K1976}. 

This paper is organized as follows. In section \ref{section2},  we prove Proposition \ref{thm:main1} and Theorem \ref{stabilization} by introducing the one-parameter families of solutions $\{ u^k_A\}, \{ u^k\}$
as in \cite{K1976}. Section \ref{section3} is to prove Theorem \ref{oscillation}.

\setcounter{equation}{0}
\setcounter{theorem}{0}

\section{Proofs of Proposition \ref{thm:main1} and Theorem \ref{stabilization}}
\label{section2}

\noindent{\sl Proof of }Proposition \ref{thm:main1}.
Note that the function $u^k_A$ defined by $u^k_A(x,t)  = u_A(kx, k^2t)$ for $k>0$ also solves problem \eqref{heat Cauchy} with $\Omega, \sigma$ replaced by $\Omega_A, \sigma_A$, respectively. Thus it follows from the uniqueness of the solution of \eqref{heat Cauchy} that the solution $u_A$ is self-similar, namely,
\begin{equation}
\label{self-smilarity}
u_A(kx, k^2t) = u_A(x,t)
\end{equation}
for every $(x,t) \in \mathbb R^N\times(0, + \infty)$ and every $k > 0$, from which we infer that
\begin{equation}
\label{constant a}
u_A(0,t) = u_A(0,1) \  \mbox{ for every  } t > 0.
\end{equation}
Note that $0< u_A <1$ in $\mathbb R^N\times(0,+\infty)$  by the maximum  principle.

By the H\"older estimate for $u_A$ as in \cite[p. 526]{Z1977}, there exist constants $C > 0$  and $0<\theta <1$ such that
\begin{equation}
\label{Holder estimate 1}
|u_A(x, 1) - u_A(0, 1)| \le C|x|^\theta\ \mbox{ for every } x \in \mathbb R^N.
\end{equation}
It then follows from \eqref{self-smilarity} that
\begin{equation}
\label{Holder estimate 2}
|u_A(x, t) - u_A(0, t)| \le C|x|^\theta t^{-\frac \theta 2}\ \mbox{ for every } (x,t)  \in \mathbb R^N \times (0, +\infty).
\end{equation}
Combining \eqref{constant a} with this yields \eqref{stabilization for self-similar solutions}.
\qed

\medskip

\noindent
{\sl Proof of {\rm Theorem \ref{stabilization}}.}\   Let $u=u(x,t)$ be the unique bounded solution  of problem \eqref{heat Cauchy}. As in \cite{EKT2009, K1976}, we introduce the one-parameter family of functions $\{u^k\}$ by
$$
u^k(x,t)=u(kx,k^2 t)\ \mbox{ for } (x,t) \in \mathbb R^N\times (0, + \infty)  \mbox{ and  } k > 0.
$$
Let us show that $u^k$ converges to the self-similar solution  $u_A$ as $k \to \infty$ uniformly in each compact set in $\mathbb R^N\times (0,+\infty)$. By the maximum  principle we have
\begin{equation}
\label{positive values between 0 and 1}
0 < u^k(x, t) < 1\ \mbox{ for every } (x,t) \in \mathbb R^N\times (0, + \infty) \ \mbox{ and every } k > 0.
\end{equation}
Moreover,  each $u^k$ solves the problem \eqref{heat Cauchy} with $\Omega, \sigma$ replaced by $\Omega^k, \sigma^k$, respectively, where
$$
\Omega^k=\{ x\in \mathbb R^N : kx\in\Omega\} \quad \mbox{and} \quad \sigma^k =
\begin{cases}
\sigma_+ \quad&\mbox{in } \Omega^k, \\
\sigma_- \quad &\mbox{in } \mathbb R^N \setminus \Omega^k.
\end{cases}
$$

As in \cite[Lemmas 1 and 2]{K1976}, we have the following two lemmas which come from \eqref{positive values between 0 and 1} together with the H\"older and energy estimates for solutions of parabolic equations of second order with discontinuous coefficients (see,  for example, \cite[Chapter III]{LSU1968}).

\begin{lemma}
\label{equicontinuity}
For every $\tau > 0$, there exist  two constants $C_1 > 0$  and $0< \theta <1$  such that
$$
|u^k(x,t)-u^k(y,s)|\le C_1(|x-y|^\theta+|t-s|^{\frac \theta2})
$$
for every $(x,t), (y,s)\in \mathbb R^N\times [\tau,+\infty)$ and for every $k>0$.
\end{lemma}

\begin{lemma}
\label{energy estimate} For every $T > 0$ and $\rho > 0$ there exists a constant $C_2 > 0$ such that
$$
\int_0^T\!\!\int_{B_\rho(0)}|\nabla u^k|^2dxdt\le C_2
$$
for every $k> 0$. Here, $B_\rho(0)=\{x\in\mathbb R^N : |x|<\rho\}$.
\end{lemma}

Since the assumption \eqref{cone} means that $\partial\Omega$ lies between the two conical surfaces, $\partial\Omega_A$ and its translate $ \partial(\Omega_A - h p)$, this geometric condition, combined with  the homothety of $\Omega^k$ and $\Omega$ with ratio $k$,
implies directly that $\mathcal X_{\Omega^k}(x), \sigma^k(x)$ converge to $\mathcal X_{\Omega_A}(x), \sigma_A(x)$, respectively,  as $k \to \infty$,  for  every $x\in\mathbb R^N$.  Therefore, by using the diagonal process with the aid of the estimates \eqref{positive values between 0 and 1}, Lemma \ref{equicontinuity} and Lemma \ref{energy estimate},  we infer from the uniqueness of the solution of  problem \eqref{heat Cauchy} that $u^k$ converges to $u_A$ as $k \to \infty$ uniformly  in each compact set in $\mathbb R^N\times (0,+\infty)$.  Hence in particular   $u(0,t)$ converges to $u_A(0,1)$ as $t \to \infty$. Moreover it follow from Lemma \ref{equicontinuity} that
$$
|u(x,t)-u(0,t)|\le C_1|x|^\theta t^{-\frac \theta 2}\ \mbox{ for every }(x,t) \in \mathbb R^N\times (0,+\infty),
$$
which yields the desired conclusion \eqref{limita}.  \qed

\setcounter{equation}{0}
\setcounter{theorem}{0}

\section{Proof of  Theorem \ref{oscillation}}
\label{section3}

We utilize the Gaussian bounds for the fundamental solutions of diffusion equations due to
Aronson \cite[Theorem 1, p. 891]{A1967} (see also \cite[p. 328]{FS1986}). Let $g = g(x,\xi,t)$ be the fundamental solution of $u_t=\mbox{ div}(\sigma\nabla u)$. Then there exist two positive constants $\lambda < \Lambda$ such that
\begin{equation}
\label{Gaussian bounds}
\lambda t^{-\frac N2}e^{-\frac{|x-\xi|^2}{\lambda t}}\le g(x,\xi,t) \le \Lambda t^{-\frac N2}e^{-\frac{|x-\xi|^2}{\Lambda t}}
\end{equation}
for all $(x,t), (\xi,t) \in \mathbb R^N\times(0,+\infty)$, where  the constants $\lambda, \Lambda$ depend only on $N$ and the bounds $m, M$ of $\sigma$.

To construct the domain $\Omega$ with the desired property, we first choose two domains $A, B$ in $ S^{N-1}$ such that $A\subset B \subsetneq S^{N-1}$ and
\begin{equation}\label{the two domains in the unit sphere}
\mathcal H^{N-1}(B)\lambda^{\frac {N+2}2} > \mathcal H^{N-1}(A)\Lambda^{\frac {N+2}2},
\end{equation}
where $\mathcal H^{N-1}$ is the standard $(N-1)$-dimensional Hausdorff measure. Let
\begin{equation}\label{spherical area}
\alpha := \mathcal H^{N-1}(A) \quad\mbox{and}\quad \beta :=\mathcal H^{N-1}(B)
\end{equation}
for simplicity of expression. According to \eqref{the two domains in the unit sphere}, we may choose a small number $\varepsilon$ with $0 < \varepsilon < 1$ so that
\begin{equation}
\label{a key inequality}
\left[(1-\varepsilon) \beta+\varepsilon\alpha\right]\lambda^{\frac {N+2}2} > \left[(1-\varepsilon)\alpha+\varepsilon\beta\right]\Lambda^{\frac {N+2}2}.
\end{equation}
Let $\delta, R$ be two numbers such that $0<\delta < 1 < R$ and
\begin{equation}
\label{essence}
\int_\delta^{\delta R}e^{-s^2}s^{N-1}ds = (1-\varepsilon) \int_0^\infty e^{-s^2}s^{N-1}ds.
\end{equation}
We then define a sequence of numbers $\{r_n\}$ by $r_0=0$ and
\begin{equation}\label{rndef}
r_n=\delta R^{n-1}, \quad n \in \mathbb N,
\end{equation}
and a sequence of  sets $\{E_n\}$ in $\mathbb R^N$ by
\begin{align*}
E_{2(k-1)} &= \{ x \in \mathbb R^N : x = r\omega,\ r_{2(k-1)}\le r \le r_{2k-1},\ \omega \in A \},\\
E_{2k-1} &= \{ x \in \mathbb R^N : x = r\omega,\ r_{2k-1}\le r \le r_{2k},\ \omega \in B \}
\end{align*}
for $k \in \mathbb N$. At last, we define the domain $\Omega$ to be the interior of the set $\bigcup\limits_{n=0}^\infty E_n$.

Since the initial condition of the problem \eqref{heat Cauchy} is the characteristic function of $\Omega$, the solution $u(x,t)$ is given by
$$
u(x,t)= \int_{\Omega} g(x,\xi, t) d\xi.
$$
It follows from \eqref{Gaussian bounds} that for every $t > 0$
\begin{equation}\label{estimates from above and below}
\lambda t^{-\frac N2}\int_\Omega e^{-\frac{|\xi|^2}{\lambda t}}d\xi\le u(0,t) \le   \Lambda t^{-\frac N2}\int_\Omega e^{-\frac{|\xi|^2}{\Lambda t}} d\xi.
\end{equation}
Let us calculate the both sides of \eqref{estimates from above and below}. Using the notation \eqref{spherical area}, we have
\begin{align}
\lambda t^{-\frac N2}\int_\Omega e^{-\frac{|\xi|^2}{\lambda t}}d\xi
&=\lambda t^{-\frac N2} \sum_{k=1}^\infty \left[ \int_{E_{2(k-1)}} e^{-\frac{|\xi|^2}{\lambda t}}d\xi + \int_{E_{2k-1}} e^{-\frac{|\xi|^2}{\lambda t}}d\xi \right] \nonumber
\\
& =\lambda t^{-\frac N2} \sum_{k=1}^\infty \left[ \alpha\int_{r_{2(k-1)}}^{r_{2k-1}} e^{-\frac{r^2}{\lambda t}}r^{N-1}dr+ \beta \int_{r_{2k-1}}^{r_{2k}}e^{-\frac{r^2}{\lambda t}}r^{N-1}dr \right] \nonumber
\\
&=\lambda^{\frac {N+2}2} \sum_{k=1}^\infty \left[ \alpha \int_{\frac {r_{2(k-1)}}{\sqrt{\lambda t}}}^{\frac {r_{2k-1}}{\sqrt{\lambda t}}} e^{-s^2}s^{N-1}ds+ \beta \int_{\frac {r_{2k-1}}{\sqrt{\lambda t}}}^{\frac {r_{2k}}{\sqrt{\lambda t}}}e^{-s^2}s^{N-1}ds \right]. \qquad \label{series expression left}
\end{align}
Replacing $\lambda$ by $\Lambda$ yields that
\begin{equation}
\label{series expression right}
\Lambda t^{-\frac N2}\int_\Omega e^{-\frac{|\xi|^2}{\Lambda t}}d\xi=\Lambda^{\frac {N+2}2} \sum_{k=1}^\infty \left[ \alpha \int_{\frac {r_{2(k-1)}}{\sqrt{\Lambda t}}}^{\frac {r_{2k-1}}{\sqrt{\Lambda t}}} e^{-s^2}s^{N-1}ds+ \beta \int_{\frac {r_{2k-1}}{\sqrt{\Lambda t}}}^{\frac {r_{2k}}{\sqrt{\Lambda t}}}e^{-s^2}s^{N-1}ds \right].
\end{equation}

Let us consider the sequence of times $\{t_n\}$ defined by
$$
 t_n=\frac 1\lambda R^{4(n-1)}, \quad n\in\mathbb N.
$$
According to \eqref{estimates from above and below} and \eqref{series expression left}, we have
\begin{equation}\label{sum of sequences}
u(0,t_n)\ge \lambda^{\frac {N+2}2} \left[ \alpha \int_0^\infty e^{-s^2}s^{N-1}ds + (\beta-\alpha)  \sum_{k=1}^\infty \int_{\frac {r_{2k-1}}{\sqrt{\lambda t_n}}}^{\frac {r_{2k}}{\sqrt{\lambda t_n}}} e^{-s^2}s^{N-1}ds \right].
\end{equation}
Since $\beta > \alpha$, we have
$$
u(0,t_n)\ge \lambda^{\frac {N+2}2} \left[ \alpha \int_0^\infty e^{-s^2}s^{N-1}ds + (\beta-\alpha)
\int_{\frac {r_{2n-1}}{\sqrt{\lambda t_n}}}^{\frac {r_{2n}}{\sqrt{\lambda t_n}}}e^{-s^2}s^{N-1}ds \right].
$$
By the definition \eqref{rndef} of $r_l$, $\frac {r_{2n-1}}{\sqrt{\lambda t_n}}=\delta$ and $\frac {r_{2n}}{\sqrt{\lambda t_n}}=\delta R$, and hence we have
$$
\int_{\frac {r_{2n-1}}{\sqrt{\lambda t_n}}}^{\frac {r_{2n}}{\sqrt{\lambda t_n}}}e^{-s^2}s^{N-1}ds
= \int_{\delta}^{\delta R} e^{-s^2}s^{N-1}ds = (1-\varepsilon) \int_0^\infty e^{-s^2}s^{N-1}ds,
$$
where the second identity follows from \eqref{essence}.
Thus we have
$$
u(0,t_n)\ge \lambda^{\frac {N+2}2}\left( (1-\varepsilon)\beta + \varepsilon \alpha\right) \int_0^\infty e^{-s^2}s^{N-1}ds
$$
for every $n \in\mathbb N$, and hence
\begin{equation}
\label{estimate from below}
\limsup_{t \to \infty}u(0,t)\ge  \lambda^{\frac {N+2}2}\left( (1-\varepsilon)\beta + \varepsilon \alpha\right) \int_0^\infty e^{-s^2}s^{N-1}ds.
\end{equation}

We now consider  the sequence of times $\{T_n\}$ defined by
 $$
 T_n=\frac 1\Lambda R^{2(2n-3)}, \quad n=2,3, \ldots.
 $$
In the same way as above, one can show  using \eqref{series expression right} that
$$
u(0,T_n)\le \Lambda^{\frac {N+2}2}\left( (1-\varepsilon) \alpha+\varepsilon\beta\right) \int_0^\infty e^{-s^2}s^{N-1}ds,
$$
and hence
\begin{equation}
\label{estimate from above}
\liminf_{t \to \infty}u(0,t)\le  \Lambda^{\frac {N+2}2}\left(  (1-\varepsilon) \alpha+\varepsilon\beta\right) \int_0^\infty e^{-s^2}s^{N-1}ds.
\end{equation}
Therefore, it follows from \eqref{a key inequality}, \eqref{estimate from below} and \eqref{estimate from above} that
\begin{equation}
\label{different}
\liminf_{t \to \infty}u(0,t) < \limsup_{t \to \infty}u(0,t).
\end{equation}

Since $\alpha < \beta$, we have from \eqref{estimates from above and below}, \eqref{series expression left} and \eqref{sum of sequences} that for every $t>0$
\begin{equation}
\label{estimate from below for all time}
u(0,t)\ge  \lambda^{\frac {N+2}2} \alpha \int_0^\infty e^{-s^2}s^{N-1}ds,
\end{equation}
which shows that
\begin{equation}
\label{strict positivity}
\liminf_{t \to \infty}u(0,t)  > 0.
\end{equation}
Observe that $1-u$ solves problem \eqref{heat Cauchy} where $\Omega$ is replaced by $\Omega^c=\mathbb R^N\setminus\Omega$, and the fact that $\Omega_{B^c}\subset\Omega^c$ with $B^c=S^{N-1}\setminus B$.
Here we used the notation \eqref{definition of cone} for $\Omega_{B^c}$.
Hence,  by the same argument as that for \eqref{estimate from below for all time},  we infer that for every $t > 0$
$$
1-u(0,t)\ge  \lambda^{\frac {N+2}2} \mathcal H^{N-1}(B^c) \int_0^\infty e^{-s^2}s^{N-1}ds,
$$
which yields that
\begin{equation}
\label{strictly less than 1}
\limsup_{t \to \infty}u(0,t)  < 1.
\end{equation}
This completes the proof.  \qed


\end{document}